\documentclass[10pt,bezier]{article}
\usepackage{amsmath,graphicx,amssymb,amsfonts,color}

\font\bg=cmbx10 scaled\magstep1

\font\small=cmr8

\newtheorem{newlemma}{{\bf Lemma}}

\newenvironment{lema}{\begin{newlemma}{\hspace{-0.5
em}{\bf.}}}{\end{newlemma}}

\newtheorem{newteorem}{{\bf Theorem}}

\newenvironment{teorem}{\begin{newteorem}{\hspace{-0.5
em}{\bf.}}}{\end{newteorem}}

\newtheorem{newkorolari}{{\bf Corollary}}

\newenvironment{korolari}{\begin{newkorolari}{\hspace{-0.5
em}{\bf.}}}{\end{newkorolari}}

\newtheorem{newdefine}{{\bf Definition}}

\newenvironment{define}{\begin{newdefine}{\hspace{-0.5
em}{\bf.}}}{\end{newdefine}}

\newtheorem{newquestion}{{\bf Question}}

\newtheorem{newkonjek}{{\bf Conjecture}}

\newtheorem{newexample}{{\bf Example}}

\begin{document}
\tolerance=10000
\baselineskip18truept
\newbox\thebox
\global\setbox\thebox=\vbox to 0.2truecm{\hsize
0.15truecm\noindent\hfill}
\def\boxit#1{\vbox{\hrule\hbox{\vrule\kern0pt
     \vbox{\kern0pt#1\kern0pt}\kern0pt\vrule}\hrule}}
\def\qed{\lower0.1cm\hbox{\noindent \boxit{\copy\thebox}}\bigskip}
\def\ss{\smallskip}
\def\ms{\medskip}
\def\bs{\bigskip}
\def\c{\centerline}
\def\nt{\noindent}
\def\ul{\underline}
\def\ol{\overline}
\def\lc{\lceil}
\def\rc{\rceil}
\def\lf{\lfloor}
\def\rf{\rfloor}
\def\ov{\over}
\def\t{\tau}
\def\th{\theta}
\def\k{\kappa}
\def\l{\lambda}
\def\L{\Lambda}
\def\g{\gamma}
\def\d{\delta}
\def\D{\Delta}
\def\e{\epsilon}
\def\lg{\langle}
\def\rg{\rangle}
\def\p{\prime}
\def\sg{\sigma}
\def\ch{\choose}

\newcommand{\ben}{\begin{enumerate}}
\newcommand{\een}{\end{enumerate}}
\newcommand{\bit}{\begin{itemize}}
\newcommand{\eit}{\end{itemize}}
\newcommand{\bea}{\begin{eqnarray*}}
\newcommand{\eea}{\end{eqnarray*}}
\newcommand{\bear}{\begin{eqnarray}}
\newcommand{\eear}{\end{eqnarray}}

\centerline{\Large  On the number of outer connected dominating }
 \vspace{.3cm}

\centerline {\Large  sets of  graphs}
\bigskip

\bigskip

\baselineskip12truept\vskip-9truept \centerline{  S.
Alikhani$^{a}$
, M. H.
Akhbari$^{b}$, C. Eslahchi$^c$  and R. Hasni$^d$ }
\baselineskip20truept \baselineskip20truept \centerline{\it
$^{a}$Department of Mathematics, Yazd University} \vskip-8truept
\centerline{\it 89195-741, Yazd, Iran}

 \centerline{\it $^{b}$School of Mathematical Sciences}
\vskip-8truept \centerline{\it Universiti Sains Malaysia, 11800
USM, Penang, Malaysia} 
\centerline{\it
$^c$Department of Mathematics, Shahid Beheshti
University}\vskip-8truept \centerline{\it G. C. Tehran, Iran }

\centerline{\it $^d$Department of Mathematics, Faculty of Science
and Technology,}  \vskip-8truept \centerline{\it Universiti
Malaysia Terengganu, 21030  Kuala Terengganu, Malaysia}

\nt\rule{12cm}{0.1mm}

\thispagestyle{empty}
\nt{\bg ABSTRACT}
 \medskip

\noindent{ Let $G=(V,E)$ be a simple graph.  A set $S\subseteq
V(G)$ is called an outer-connected dominating set (or ocd-set) of
$G$, if $S$ is a dominating set of $G$ and either $S=V(G)$ or
$V\backslash S$ is a connected graph. In this paper we introduce a
polynomial which its coefficients are the number of ocd-sets of
$G$.  We obtain some properties  of this polynomial and its
coefficients. Also we compute this polynomial for some specific
graphs. }

\bs

\nt{\it Keywords:} {\small Ocd polynomial;  outer connected
dominating set; coefficient}

\nt{\it Mathematics subject classification:} {\small 05C69, 11B83}

\baselineskip20truept

\setcounter{page}{1}
\section{Introduction}
\nt Let $G=(V,E)$ be a simple graph. A set $S\subseteq V$ is a
dominating set if  every vertex in $V\backslash S$ is adjacent to
at least one vertex in $S$.
 The domination number $\gamma(G)$ is the minimum cardinality of a dominating set in $G$.

\nt A set $S\subseteq V(G)$ is called an outer-connected
dominating set (or ocd-set) of $G$, if $S$ is a dominating set of
$G$ and either $S=V(G)$ or $V \backslash S$ is a connected graph.
The outer-connected domination number $\widetilde{\gamma}_{c}$ of
$G$ is the minimum cardinality of a outer-connected dominating set
of $G$. (\cite{ajc,utilit})


\ms

\nt The {\it corona} of two graphs $G_1$ and $G_2$, as defined by
Frucht and Harary in~\cite{harary},
 is the graph
$G=G_1 \circ G_2$ formed from one copy of $G_1$ and $|V(G_1)|$
copies of $G_2$, where the ith vertex of $G_1$ is adjacent to
every vertex in the ith copy of $G_2$. The corona $G\circ K_1$, in
particular, is the graph constructed from a copy of $G$,
 where for each vertex $v\in V(G)$, a new vertex $v'$ and a pendant edge $vv'$ are added.
 The {\it join} of two graphs $G_1$ and $G_2$, denoted by $G_1\vee G_2$,
 is a graph with vertex set  $V(G_1)\cup V(G_2)$
 and edge set $E(G_1)\cup E(G_2)\cup \{uv| u\in V(G_1)$ and $v\in V(G_2)\}$.

\nt As usual we denote the complete graph, the cycle, and the path
of order $n$ by $K_n$, $C_n$ and $P_n$, respectively. Also we call
$K_{1,n}$ a star of order $n+1$.

\ms

\nt In the next section, we introduce a polynomial for a graph $G$
which its coefficients are the number of outer-connected
dominating sets of $G$  and call it the outer connected domination
polynomial. We obtain some of its properties. In Section 3, we
study this polynomial for some certain graphs. In the last section
we define $\widetilde{D}$-equivalence classes of graphs and
investigate it for some graphs.

\bs

\section{Introduction to outer-connected domination polynomial}

\nt Similar to domination polynomial of a graph (\cite{euro,ars}),
in this section we state the definition of outer-connected
domination (ocd)  polynomial and obtain some of its properties.

\begin{define}\label{Definition2}
 Let $\widetilde{{\cal D}}(G,i)$ be the family of outer-connected dominating sets of a graph $G$ with cardinality $i$ and let
$\widetilde{d}(G,i)=|{\widetilde{\cal D}}(G,i)|$. Then the outer
connectivity domination polynomial $\widetilde{D}(G,x)$ of $G$ is
defined as
\begin{center}
$\widetilde{D}(G,x)=\displaystyle\sum_{i=\widetilde{\gamma}_{c}(G)}^{|V(G)|}
\widetilde{d}(G,i) x^{i}$,
\end{center}
where $\widetilde{\gamma}_{c}(G)$ is the ocd number of $G$.
\end{define}

\nt The path $P_4$ on $4$ vertices, for example, has one
outer-connected  dominating set of cardinality $4$, four
outer-connected
 dominating sets of cardinalities  $3$ and one outer-connected
dominating set of cardinality $2$; its  ocd polynomial is then
$\widetilde{D}(P_4,x) = x^4+4x^3+x^2$. As another example, it is
easy to see that, for every $n\in \mathbb{N}$,
$\widetilde{D}(K_n,x)=(1+x)^n-1$.

\begin{teorem}\label{theorem2.2.2}{\rm(\cite{ajc})}
If a graph $G$ consists of $m$ components $G_1,\ldots,G_m$, then
$$\widetilde{\gamma}_{c}(G)=|V(G)|-\max \{|V(G_i)|-\widetilde{\gamma}_{c}(G_i):i=1,\ldots
m\}.$$
\end{teorem}

\begin{teorem}\label{theorem2.2.2}
If a graph $G$ consists of $m$ components $G_1,\ldots,G_m$, then
$\widetilde{D}(G,x)=(\sum_{i=1}^m x^{b_i}
\widetilde{D}(G_i,x))-(m-1)x^n$, where $n=|V(G)|$ and
$b_i=n-|V(G_i)|$.
\end{teorem}

\nt{\bf Proof.} Let $A$ be an ocd set of graph $G$. Obviously we
have $A=B_i\bigcup _{j\neq i} V(G_j)$, where $B_i$ is an ocd set
of $G_i$ for $1\leq i\leq m$. Note that the coefficient of $x^n$
in $\widetilde{D}(G,x)$ is one, but in $\sum_{i=1}^m x^{b_i}
\widetilde{D}(G_i,x)$, we count $x^n$, $m$-times. So we should
have $-(m-1)x^n$ in the right of equality.\quad\qed

\nt As a consequence of Theorem~\ref{theorem2.2.2},
 we have the following corollary for the empty graphs:

\begin{korolari}\label{corollary2.2.3}
Let $\overline{K}_n$ be the empty graph with $n$ vertices. Then
$\widetilde{D}(\overline{K}_n,x)=x^n$.
\end{korolari}

\nt{\bf Proof.} Since $\widetilde{D}(\overline{K}_1,x)=x$, we have
the result by Theorem~\ref{theorem2.2.2}.\quad\qed

\begin{korolari}
Let $G$ be a graph with $r$ isolated vertices $v_1,\ldots,v_r$.
Then
$$\widetilde{D}(G,x)=x^r \widetilde{D}(G-\{v_1,\ldots,v_r\},x).$$
\end{korolari}

\ms

\nt The following theorem is an easy consequence of the definition
of the ocd polynomial.

\begin{teorem}\label{theorem3}
Let $G$ be a graph with $|V(G)|=n$. Then
\begin{enumerate}
\item[(i)] If $G$ is connected, then $\widetilde{d}(G,n)=1$ and
$\widetilde{d}(G,n-1)=n$, \item[(ii)] $\widetilde{d}(G,i)=0$ if
and only if $i<\widetilde{\gamma}(G)$ or $i>n$. \item[(iii)]
$\widetilde{D}(G,x)$ has no constant term. \item[(iv)]
$\widetilde{D}(G,x)$ is a strictly increasing function in
$[1,\infty)$. \item[(v)] Let $G$ be a graph and $H$ be any induced
subgraph of $G$. Then $deg(\widetilde{D}(G,x))\geq
deg(\widetilde{D}(H,x))$.
 \item[(vi)] Zero is a root of $\widetilde{D}(G,x)$,
with multiplicity $\widetilde{\gamma}_{c}(G)$.
 \item[(vii)] For
any spanning subgraph $H$ of $G$, $\widetilde{\gamma}_c(G)\leq
\widetilde{\gamma}_c(H)$.
\end{enumerate}
\end{teorem}

\bs

\section{OCD polynomial of certain graphs}

\nt In this section, we obtain some formulas  for computing the
ocd polynomial of some certain graphs. We recall the following
theorem.

\ms \nt

\begin{teorem}\label{theorem5'}{\rm(\cite{ajc})}
\begin{enumerate}
\item[(i)] $ \widetilde{\gamma}_c(P_n)=\left\{
\begin{array}{lr}
{\displaystyle n-1};&
\quad\mbox{$n=2,3$,}\\[15pt]
{\displaystyle n-2};& \quad\mbox{ $n\geq 4$;}
\end{array}
\right. $

\item[(ii)] For every $n\geq 3$, $\widetilde{\gamma}_c(C_n)=n-2$.

\item[(iii)] If $t\geq 2$ and $n_1\leq n_2\leq \ldots \leq n_t$,
then

$ \widetilde{\gamma}_c(K_{n_1,\ldots,n_t})=\left\{
\begin{array}{lr}
{\displaystyle n_2};&
\quad\mbox{if $t=2$ and  $n_1=1$,}\\[15pt]
{\displaystyle 1};& \quad\mbox{if $t\geq 3$  and $n_1=1$.}\\[15pt]
{\displaystyle 2}; & \quad\mbox{if $t\geq 2$ and $n_1>1$.}
\end{array}
\right. $

\item[(iv)] If $G$ is a connected graph of order $n\geq 2$,  then
$\widetilde{\gamma}_c(G)=n-1$ if and only if $G$ is a star.

\end{enumerate}
\end{teorem}

\begin{teorem}\label{theorem5}
\begin{enumerate}
\item[(i)] For every $n\geq 4$,
$\widetilde{D}(P_n,x)=x^{n-2}(x^2+nx+n-3)$.

\item[(ii)] For every $n\geq 3$,
$\widetilde{D}(C_n,x)=x^{n-2}(x^2+nx+n)$.

\item[(iii)] For every $n\in \mathbb{N}$,
$\widetilde{D}(K_{1,n},x)=x^n(n+x+1)$.

\end{enumerate}
\end{teorem}
\nt{\bf Proof.}
\begin{enumerate}
\item[(i)] It is obvious that $\widetilde{d}(P_n,n)=1$, and
$\widetilde{d}(P_n,n-1)=n$. Also one can see that
$d(P_n,n-2)=n-3$. Since by Theorems~\ref{theorem3} and
\ref{theorem5'}(i), $\widetilde{d}(P_n,k)=0$ for $k\leq n-3$,  we
have the result.

\item[(ii)] It is obvious that $\widetilde{d}(C_n,n)=1$, and
$\widetilde{d}(C_n,n-1)=n$. Every outer-connected dominating set
of $C_n$ which labeled by $\{1,\ldots,n\}$ is one of  the form of
set
$\Big\{\{1,\ldots,n-2\},\{2,\ldots,n-1\},\ldots,\{n,1,\ldots,n-2\}\Big\}$
so $d(C_n,n-2)=n$. Since by Theorems~\ref{theorem3} and
\ref{theorem5'}(ii), $\widetilde{d}(C_n,k)=0$ for $k\leq n-3$,
Therefore we have the result.

\item[(iii)] It follows from Theorem \ref{theorem3}(ii) and
\ref{theorem5'}(iv).\quad\qed

\end{enumerate}

\nt Here, we provide a formula for the  ocd polynomial of the join
of two graphs. First we state the following lemmas which are about
the outer-connected domination number of join of two graphs.

\begin{lema}\label{lemma10'}
Let $G_1$ be a graph of order $n\geq 2$. If $G$ has a vertex $v$
of degree $n-1$, then for every graph $G_2$,
$\widetilde{\gamma}_c(G_1 \vee G_2)=1$.
\end{lema}

\nt{\bf Proof.} It is easy to see that $\{v\}$ is an ocd set of
$G$. Therefore we have the result.\quad\qed

\begin{lema}\label{lemma10} Let $G_1$ be any graph which does not have a vertex with degree $|V(G_1)|-1$, and $G_2$ be an arbitrary graph. Then
$\widetilde{\gamma}_c(G_1 \vee G_2)=2.$
\end{lema}
\nt{\bf Proof.} Suppose that $v\in V(G_1)$ and $w\in V(G_2)$. It
is easy to see that $\{v,w\}$ is an outer-connected dominating set
of $G_1 \vee G_2$. So  $\widetilde{\gamma}_c(G_1 \vee G_2)\leq 2$.
Now suppose that  there exist a vertex $u\in V(G_i)$ $(i=1,2)$
such that $\{u\}$ is an ocd of $G_1 \vee G_2$. In this case we
have $deg u=|V(G_i)|-1$ and this a contradiction. Therefore
$\widetilde{\gamma}_c(G_1 \vee G_2)=2.$\quad\qed

\begin{teorem}\label{theorem2.2.7}
Let $G_1$ and $G_2$ be two connected  graphs of order  $m$ and
$n$, respectively. If $n_i$ $(i=1,2)$ is the number of vertices of
$G_i$ $(i=1,2)$ with degree $|V(G_i)|-1$, then
\[
\widetilde{D}(G_1\vee G_2,x)=(1+x)^{m+n}+(n_1+n_2-m-n)x-1.
\]
\end{teorem}

\nt{\bf Proof.} Let $i$ be a natural number  $1\leq i \leq m+n$.
We want to determine $\widetilde{d}(G_1\vee G_2,i)$. If $i_1$ and
$i_2$ are two natural numbers such that $i_1+i_2=i$,
 then  for every $D_1\subseteq  V(G_1)$ and $D_2\subseteq V(G_2)$,
such that $|D_j|=i_j$, $j=1,2$, $D_1\cup D_2$ is an ocd set of
$G_1\vee G_2$.  On the other hand by Lemmas \ref{lemma10'} and
\ref{lemma10}  one vertex $x$ of $G_i$ $(i=1,2)$ is an ocd of
$G_1\vee G_2$, if $x$ is adjacent to all vertices of $G_i$.
Therefore we have the result.\quad\qed

\nt As a corollary, we have the following formula for the ocd
polynomial of 
the wheel $W_n$.

\begin{korolari}\label{corollary2.2.8}
  If $n\geq 4$,
then
$$\widetilde{D}(W_n,x)=(1+x)^{n}-(n+1)x-1.$$

\end{korolari}

\nt {\bf Proof.}
Since for every $n\geq 4$, $W_n=C_{n-1}\vee K_1$, we have the
result by Theorem \ref{theorem2.2.7}.\quad\qed


\nt {\bf Acknowledgements.} The first author would like to thank
the School of Mathematical Sciences at Universiti Sains Malaysia
(USM) for partial support and hospitality during his visit. Also
Authors wish to thank the referee for his/her valuable comments.

\end{document}